\documentclass[10pt]{amsart}
\usepackage{a4wide}
\usepackage{amsfonts}
\usepackage{psfrag}
\usepackage{bbm}

\usepackage{amsmath}
\usepackage{amssymb}
\usepackage{amsthm}
\usepackage{graphicx}

\newtheorem {thm}{Theorem}
\newtheorem {rem}{Remark}

\newtheorem {cor}{Corollary}
 
\newcommand{\E}{\mathbb{E}}

\renewcommand{\P}{\mathbb{P}}
\numberwithin{equation}{section}
\numberwithin{equation}{section}
\numberwithin{equation}{section}
\newcommand{\Z}{\mathbb{Z}}

\begin{document}

\title{An Application of Renewal Theorems to Exponential Moments of Local Times}
\author{Leif D\"oring}

\address{Institut f\"ur Mathematik, Technische Universit\"at Berlin, Stra\ss e des 17.~Juni 136, 10623 Berlin, Germany}
\thanks{The first author was supported by the EPSRC grant EP/E010989/1.}
\email{leif.doering@googemail.com}
\author{Mladen Savov}
\address{New College, University of Oxford, Holywell Street, Oxford OX1 3BN, United Kingdom}
\email{savov@stats.ox.ac.uk}
\date{}

\maketitle
\maketitle

\begin{abstract}
In this note we explain two transitions known for moment generating functions of local times by means of properties of the renewal measure of a related renewal equation. The arguments simplify and strengthen results on the asymptotic behavior in the literature.
\end{abstract}

	\section{Introduction and Results}
	Suppose $X=(X_t)$ is a time-homogeneous continuous time Markov process on a countable set $S$ with transition probabilities $p_t(i,j)=\P[X_t=j\,|X_0=i]$ for $i,j\in S$. We fix some arbitrary $i\in S$ and denote by $L^i_t$ the time $(X_t)$ spends at $i$ until time $t$:
	\begin{align*}
		L_t^i=\int_0^t\delta_i(X_s)\,ds.
	\end{align*}
	 A quantity that has been studied in different contexts is the moment generating function $\E^i\big[ e^{\gamma L^i_t}\big]$, where $X_0=i$ and $\gamma$ is a positive real number. \\
	To explain our interest in $\E^i\big[ e^{\gamma L^i_t}\big]$ let us have a brief look at the parabolic Anderson model with Brownian potential, i.e.
	\begin{align}\label{pam}
		du_t(i)=\Delta u_t(i)\,dt+\gamma u_t(i)dB_t(i)
	\end{align}
	with homogeneous initial conditions $u_0\equiv 1$. Here, $i\in \Z^d$, $\Delta$ denotes the discrete Laplacian $\Delta f(i)=\sum_{|i-j|=1}1/(2d)(f(j)-f(i))$, and $\{B(i)\}_{i\in \Z^d}$ is a family of independent Brownian motions. It is known (see for instance Theorem II.3.2 of \cite{CM94}) that the moments of $u_t(i)$ solve discrete-space heat equations with one-point potentials. In particular, $\E[u_t(i)u_t(j)]$ solves
	\begin{align}\label{eq}
		\frac{d}{dt}w(t,i,j)=\Delta w(t,i,j)+\gamma\delta_0(i-j)w(t,i,j)
	\end{align}
	with homogeneous initial conditions. The discrete Laplacian acts on both spatial variables $i$ and $j$ seperately. Applying the Feynman-Kac formula one reveals that 
	\begin{align*}
		w(t,i,j)=\E^{i,j}\big[e^{\gamma \int_0^t\delta_0(X^1_s-X^2_s)\,ds}\big]
	\end{align*}
	where $X^1,X^2$ are independent simple random walks. Hence, for $L_t$ corresponding to the difference walk $X^1-X^2$ (or equivalently corresponding to a simple random walk with doubled jump rate)
	\begin{align*}
		\E\big[u_t(i)^2\big]=w(t,i,i)=\E^0\big[e^{\gamma L_t^0}\big].
	\end{align*}
	The notion of weak $2$-intermittency, i.e. exponential growth of the second moment $\E[u_t(i)^2]$, now explains the interest in the study of the exponential moment of $L_t^i$ for continuous time Markov processes. 
	
	Applying the variation of constant formula to solutions of (\ref{eq}) one can guess that the following renewal equation holds for fixed $\gamma\geq 0$ and $t\geq 0$:
	\begin{align}\label{ren}
		\E^i\big[e^{\gamma L^i_t}\big]=1+\gamma\int_0^t \E^i\big[e^{\gamma L^i_{t-s}}\big]p_s(i,i)\,ds.
	\end{align}
	Indeed, expanding the exponential one can show directly the validity of Equation (\ref{ren}) for general time-homogeneous Markov processes on countable state spaces (see Lemma 3.2 of \cite{AD09}). The same equation holds for $\E^i\big[e^{L_t^j}\big]$ with $p_s(i,i)$ replaced by $p_s(i,j)$. As the 	analysis does not change we restrict ourselves to $i=j$.\\
	
	In Section III of \cite{CM94} and as well in Lemma 1.3 and Theorem 1.4 of \cite{GdH06}  analytic techniques were applied to understand the longtime behavior of solutions of (\ref{eq}) by means of spectral properties of the discrete Laplacian with one-point potential. They showed that the exponential growth rate
	\begin{align*}
		r(\gamma):=\lim_{t\rightarrow \infty}\frac{1}{t}\log \E\big[u_t(i)^2\big]=\lim_{t\rightarrow \infty}\frac{1}{t}\log w(t,i,i)=\lim_{t\rightarrow \infty}\frac{1}{t}\log \E^0\big[e^{\gamma L^0_t}\big].
	\end{align*}
	exists and obeys the following transition in $\gamma$:
	\begin{align*}
		r(\gamma)>0 \text{ if and only if } \gamma>1/G_{\infty}(i,i),
	\end{align*}
	where $G_{\infty}(i,i)$ is the Green function $\int_0^{\infty}p_s(i,i)\,ds$. This first transition in $\gamma$ can be proved analytically, as identifying $r(\gamma)$ corresponds to identifying the smallest eigenvalue of the perturbed operator $H=\Delta+\gamma\delta_0$. As multiplication with $\gamma\delta_0$ is a one-dimensional perturbation and for the discrete Laplacian explicit formulas for eigenfunctions are available, all necessary quantities can be calculated. In particular the exponential growth rate $r(\gamma)$, in the general case of our Theorem \ref{t1} represented by the Laplace transform of the transition probabilities, has been described for the simple random walk as the unique solution of
	\begin{align*}
		\frac{2}{\gamma}=\frac{1}{(2\pi)^d}\int_{S^2}\frac{1}{\Phi(s)+r(\gamma)\gamma/2}\,ds
	\end{align*}
	where $S^2$ denotes the $d$-dimensional torus and $\Phi(s)=2\sum_{i=1}^d(1-\cos(s_i))$. Compared to this expression, our Laplace transform representation is particularly useful (and easy to prove) as it immediately provides the qualitative behavior of $r(\gamma)$ as a function of $\gamma$.\\
	
	Replacing the discrete Laplacian by a generator of a finite range random walk, in \cite{DD06} the second moment of solutions of (\ref{pam}) were analyzed via a random walk representation. For more general initial conditions this leads to a renewal equation similar to (\ref{ren}). The authors analyzed their equation (3.16) directly without appealing to the renewal theorem. In precisely the same manner as we do in the proof of our Theorem \ref{t1} one can proceed in their case and strengthen the asymptotics of their Equation (3.15).\\
	
	Assuming only that $p_t(i,i)\sim ct^{-\alpha}$ for some $\alpha>0$ (by $f \sim g$ we denote strong asymptotic equivalence $\lim f/g=1$ at infinity) a second transition was revealed in Proposition 3.12 of \cite{AD09} by a Laplace transform technique combined with Tauberian theorems: at the critical point $\gamma=1/G_{\infty}$ the growth is of linear order if and only if $\alpha>2$. As for the simple random walk on $\Z^d$ the local central limit theorem implies $p_t(i,i)\sim ct^{-d/2}$, linear growth occurs for dimensions at least $5$. \\
	
	The main goal of the following is to show how the known results easily follow from different renewal theorems utilizing the fact that Equation (\ref{ren}) is a renewal equation of the type
	\begin{align}\label{re}
		Z(t)=z(t)+\int_0^tZ(t-s)U(ds)
	\end{align}	
	with $Z(t)=\E^i\big[e^{\gamma L^i_t}\big]$, initial condition $z\equiv 1$, and renewal measure $U(ds)=\gamma p_s(i,i)\,ds$. This approach is robust as there is no need to assume any properties of the underlying Markov process (neither symmetry to obtain a self-adjoint operator, nor  polynomial decay for Tauberian theorems or finite range transitions kernels for the random walk representation).

	The two transitions will now appear in terms of whether or not the renewal measure $U$
	\begin{itemize}
		\item is a probability measure,
		\item has finite mean.
	\end{itemize}
	
	In the supercritical case $\gamma>1/G_{\infty}(i,i)$ without any further consideration we obtain the strong asymptotics of $\E^i\big[e^{\gamma L^i_t}\big]$. This of course is stronger than considering the Lyapunov exponent $r(\gamma)$ that appears in \cite{CM94}, \cite{GdH06}, and \cite{AD09} as we exclude the possible existence of a subexponential factor. With further considerations this can be proved analytically but comes here for free from the renewal theorem.\medskip
	
	For the statement of the theorem we denote by $H_{\infty}(i,i)=\int_0^{\infty}sp_s(i,i)\,ds$ the expected time of hitting of two independent copies of $X$. In contrast to the Green function $G_{\infty}$ here we count the hitting time of the entire paths not only of the paths at same time. The Laplace transform in time of $p_t(i,i)$ is denoted by $\hat p(\lambda)$, $\lambda>0$, and weak asymptotic equivalence at infinity by $f \approx g$ (i.e. there are constants such that $C_1\leq \liminf f/g\leq \limsup f/g\leq C_2$).
	\begin{thm}\label{t1}
		Suppose $(X_t)$ is a time-homogeneous Markov process on $S$ started in $i$. Then for $L^i_t=\int_0^t\delta_i(X_s)\,ds$ the following holds:
		\begin{enumerate}
			\item If $\gamma>\frac{1}{G_{\infty}(i,i)}$, then $\hat p^{-1}(1/\gamma)>0$ and
				\begin{align*}
					 \E^i\big[e^{\gamma L^i_t}\big]\sim \frac{1}{\hat p^{-1}(1/\gamma)\gamma\int_0^{\infty}se^{-\hat p^{-1}(1/\gamma)s}p_s(i,i)\,ds} e^{\hat p^{-1}(1/\gamma)t}.
				\end{align*}
			\item If $\gamma=\frac{1}{G_{\infty}(i,i)}$, then 	
					\begin{align*}
						\E^i\big[e^{\gamma L^i_t}\big] \approx \frac{t}{\gamma\int_0^t\int_s^{\infty}p_r(i,i)\,drds},
					\end{align*}
					where the weak asymptotic bounds are $1$ and $2$. If moreover $$H_{\infty}(i,i)=\int_0^{\infty}sp_s(i,i)\,ds<\infty,$$ then
					\begin{align*}
						\E^i\big[e^{\gamma L^i_t}\big] \sim \frac{t}{\gamma H_{\infty}(i,i)}.
					\end{align*}
					
			\item If $0\leq\gamma<\frac{1}{G_{\infty}(i,i)}$, then 
				\begin{align*}
					\lim_{t\rightarrow \infty} \E^i\big[e^{\gamma L^i_t}\big]=\frac{1}{1-\gamma G_{\infty}(i,i)}.
				\end{align*}
			\end{enumerate}
	\end{thm}
	
	\begin{rem}
		In the general case, we only obtained weak convergence at criticality in the previous theorem with asymptotic bounds $1$ and $2$. Under the stronger assumptions $p_t(i,i)\sim ct^{-\alpha}$, in Proposition 3.12 of \cite{AD09} strong asymptotics were obtained by Tauberian theorems. The case 			of $\alpha>2$ is contained in the second part, $\alpha\leq 1$ is contained in the first part of our previous theorem and also strong asymptotics for $\alpha\in (1,2]$ can be obtained by extended renewal theorems. Here, we can directly use the infinite mean renewal Theorem 1 of \cite {AA87} to 		obtain precisely the same strong asymptotics as of Proposition 3.12 of \cite{AD09}.		
	\end{rem}

\begin{figure}
    	\vspace{-2.5cm}
	\begin{center}
                \psfrag{a}{$\gamma$}
      	 \psfrag{b}{$r(\gamma)$}
        \includegraphics[scale=0.3]{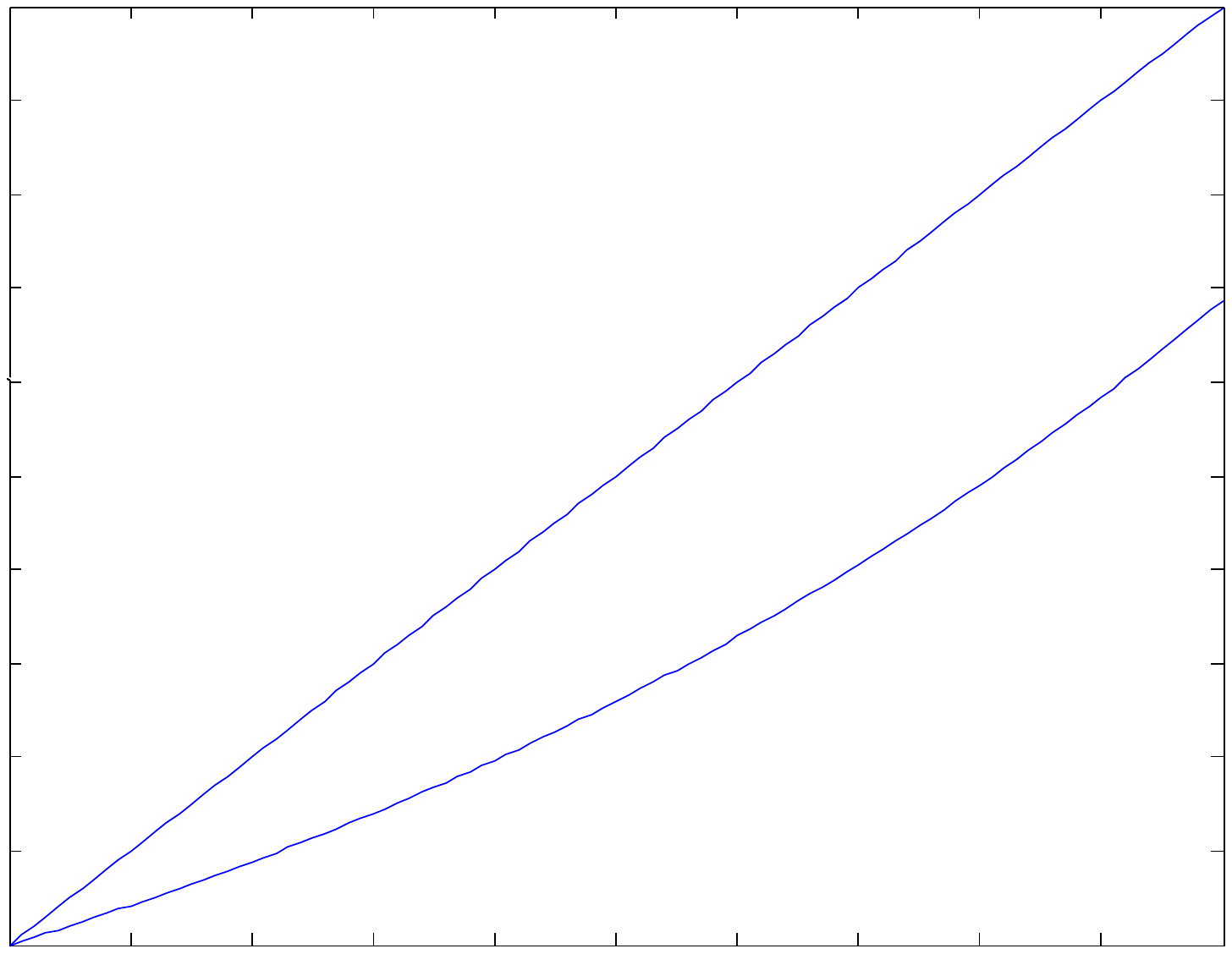}
               \psfrag{c}{$\gamma$}
       	 	\psfrag{d}{$r(\gamma)$}
        \includegraphics[scale=0.3]{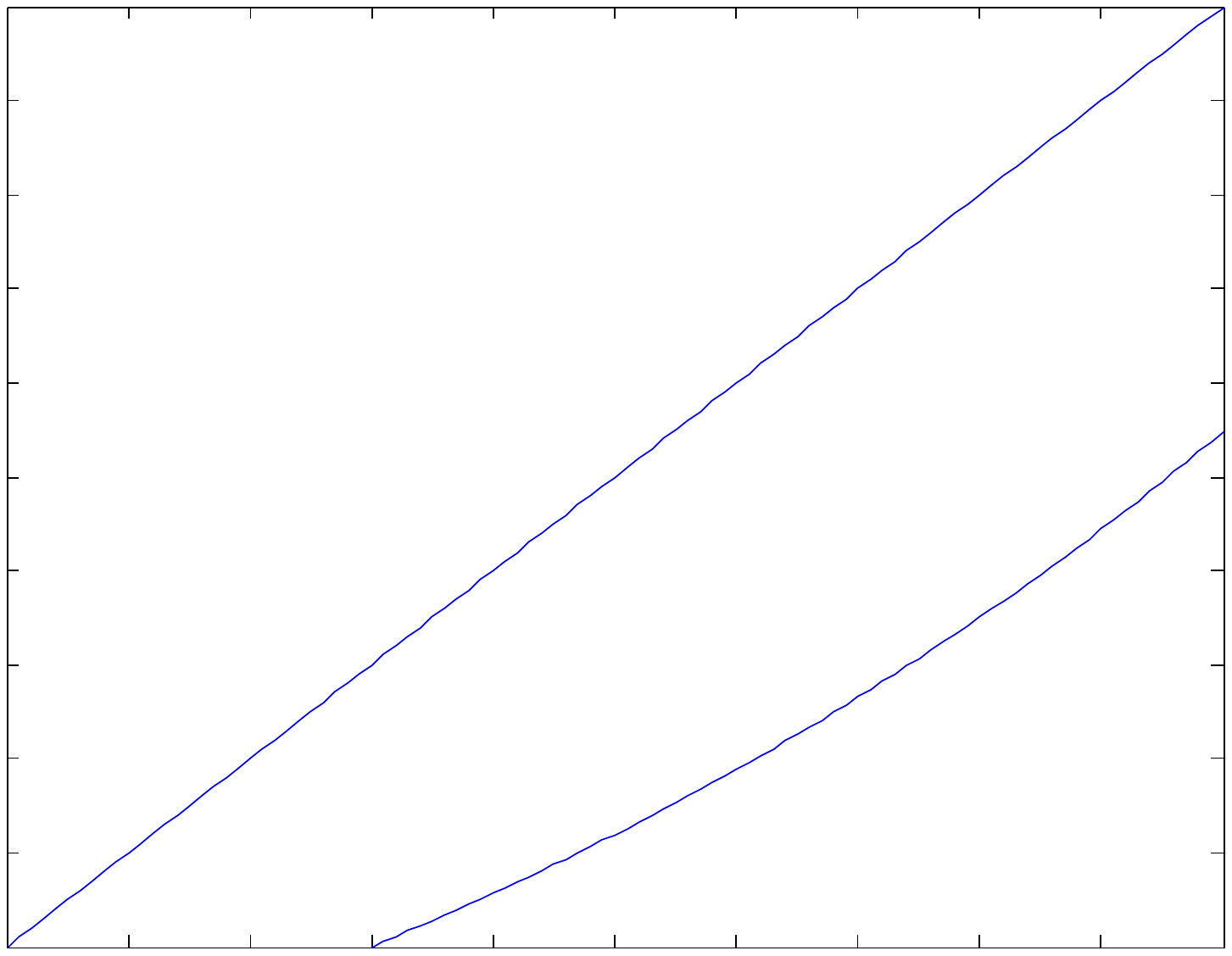}
           \end{center}
    	\vspace{-3cm}
	\caption{$\gamma\mapsto r(\gamma)$ for $G_{\infty}(i,i)=\infty$ and $G_{\infty}(i,i)<\infty$ plottet agains the identity function}
	\label{fig}
\vspace{-0.3cm}

\end{figure}

	Qualitative properties of the exponential growth rate $r(\gamma)$ have been considered for the simple random walk (see Section III of \cite{CM94}, Theorem 1.4 of \cite{GdH06}) and in the polynomially case (see Corollary 3.10 of \cite{AD09}). The representation of the growth rate in the previous 		theorem directly shows that the qualitative behavior (see Figure \ref{fig} for the qualitative behavior of $r(\gamma)$ plotted against the identity function) is valid for general Markov processes:
	\begin{cor}\label{c1}
		Suppose $(X_t)$ is a time-homogeneous Markov process on $S$ started in $i$. Then for $L^i_t=\int_0^t\delta_i(X_s)\,ds$ the following holds for $\gamma\geq 0$:
		\begin{enumerate}
			\item $r(\gamma) \geq 0$ and $r(\gamma) > 0$ if and only if $\gamma > 1/G_{\infty}(i,i)$,
			\item the function $\gamma\mapsto r(\gamma)$ is strictly convex for $\gamma>1/G_{\infty}(i,i)$,
			\item $r(\gamma)\leq\gamma$ for all $\gamma$, and $r(\gamma)/\gamma\rightarrow 1$, as $\gamma\rightarrow \infty$.
		\end{enumerate}
	\end{cor}
		
\section{Proofs}
	
	\begin{proof}[Proof of Theorem 1]
		The proof of Theorem \ref{t1} is based on the renewal equation (\ref{re}) setting $z\equiv 1$, $Z(t)=\E^i\big[e^{\gamma L_t^i}\big] $, and $U(ds)=\gamma p_s(i,i)\,ds$.\\
	1. The assumptions of the theorem directly imply that in this case $U$ is not a probability measure. Either the measure is infinite (with density bounded by $\gamma$) or it is finite with total mass strictly larger than $1$. 
		From the definition of the Laplace transform we obtain for $\lambda=\hat p^{-1}(1/\gamma)$ that
		\begin{align*}
			\bar U(ds)=\gamma e^{-\lambda s}p_s(i,i)\,ds
		\end{align*}
		is a probability measure. As by assumption $\lambda>0$, we obtain that $e^{-\lambda t}$ is directly Riemann integrable and $\bar U$ has finite mean $\gamma\int_0^{\infty}s p_s(i,i)e^{-\lambda s}\,ds$. Hence,		\begin{align*}
			e^{-\lambda t}\E^i\big[e^{\gamma L^i_t}\big]=e^{-\lambda t}+\gamma\int_0^t e^{-\lambda(t-s)}\E^i\big[e^{\gamma L^i_{t-s}}\big] e^{-\lambda s}p_s(i,i)\,ds
		\end{align*}
		is a proper renewal equation. The classical renewal theorem (see for instance page 363 of \cite{F71}) implies that
		\begin{align*}
			\lim_{t\rightarrow \infty}e^{-\lambda t}\E^i\big[e^{\gamma L^i_t}\big]= \frac{\int_0^{\infty}e^{-\lambda s}\,ds}{\gamma\int_0^{\infty}se^{-\lambda s}p_s(i,i)\,ds}=\frac{1}{\lambda\gamma\int_0^{\infty}se^{-\lambda s}p_s(i,i)\,ds}
		\end{align*}
		proving the claim.\\
	2. In the critical case $\gamma\int_0^{\infty}p_s(i,i)\,ds=1$, the measure $U$ as defined above indeed is a probability measure which does not necessarily has finite mean. Furthermore, the situation is different from the first case as now the initial condition $z\equiv 1$ is not directly Riemann 				integrable.  Iterating Equation (\ref{ren}) we obtain the representation
		\begin{align*}
			\E^i\big[e^{\gamma L^i_t}\big]=\int_0^t\sum_{n\geq 0}p_s(i,i)^{\ast n}\,ds,
		\end{align*}
		where $\ast n$ denotes $n$-fold convolutions. Note that convergence of the series is justified by boundedness of $p$. In the case of finite mean, Equation (1.2) of page 358 of \cite{F71} and the renewal theorem on page 360 now directly imply 
		\begin{align*}
			\E^i\big[e^{\gamma L^i_t}\big]\sim \frac{t}{\gamma \int_0^{\infty}sp_s(i,i)\,ds}=\frac{t}{\gamma H_{\infty}(i,i)}.
		\end{align*}
		For renewal measure $U$ with infinite mean we again use the convolution representation of $\E^i\big[e^{\gamma L^i_t}\big]$ to apply Lemma 1 of \cite{E73} showing that the denominator needs to be replaced by the truncated mean $\int_0^t\big(1-\int_0^s\gamma p_r(i,i)\,dr\big)\,ds$.\\
	3.	For $\gamma<G_{\infty}$ we may directly use the proof of Proposition 3.11 of \cite{AD09} as there no additional structure was assumed. We repeat the simple argument for completeness. Taking Laplace transform of Equation (\ref{ren}) and solving the multiplication equation in Laplace domain 		(note that under Laplace transform the convolution turns into multiplication) we obtain with $f(t)=\E^i\big[e^{\gamma L^i_t}\big]$ 
		\begin{align*}
			\hat f(\lambda)=\frac{1}{\lambda}\frac{1}{1-\gamma \hat p(\lambda)}
		\end{align*}
		for $\lambda>0$. As by assumption the second factor converges to the constant $1/(1-\gamma G_{\infty}(i,i))$ as $\lambda$ tends to zero, Karamata's Tauberian theorem (see Theorem 1.7.6 of \cite{BGT89}) implies the result. Note that as $f(t)$ is increasing, the Tauberian condition for that 			theorem is fulfilled.
	\end{proof}

\begin{proof}[Proof of Corollary 1]
	Part 1. of Theorem \ref{t1} shows that understanding $\hat p^{-1}$ suffices to understand $r(\gamma)$. This is not difficult due to the following observation: as $p$ is bounded by $1$, $\hat p(\lambda)$ is finite for all $\lambda>0$, strictly decreasing and convex with $\hat p(0) = G_{\infty}(i,i)$. Hence, $\hat p^{-1}(\lambda) = 0$ if and only if $\lambda\geq G_{\infty}(i,i)$. This implies that $\hat p^{-1}(1/\gamma)=0$ precisely for $\lambda \leq 1/G_{\infty}(i,i)$. Hence, parts 1. and 2. are proved as $r(\gamma)=\hat p^{-1}(1/\gamma)$.
	
	First note that the first part of 3. is immediate as $L_t^i\leq t$. Continuity of $p$ and $p_0(i,i)=1$ imply that  for $\epsilon>0$ there is $t_0(\epsilon)$ such that $p_t(i,i)\geq1-\epsilon$ for $t\leq t_0(\epsilon)$. Hence,
         \begin{align*}
   		\frac{1}{\gamma}&=\hat{p}(r(\gamma))=\int_0^{\infty}e^{-r(\gamma)t}p_t(i,i)\,dt
 	    	\geq (1-\epsilon)\int_0^{t_0(\epsilon)}e^{-r(\gamma)t}\,dt
	     	=(1-\epsilon)\frac{1}{r(\gamma)}\big(1-e^{-r(\gamma)t_0(\epsilon)}\big).
         \end{align*}
 	Since $r(\gamma)\rightarrow \infty$ for $\gamma\rightarrow \infty$ we obtain
	\begin{eqnarray*}
		 \liminf_{k\rightarrow \infty}\frac{r(\gamma)}{\gamma}\geq 1.
	\end{eqnarray*}
	This combined with the first part of 3. proves the second part.
\end{proof}

\section{Related Work}
	After submission of this paper the authors learned about an unpublished manuscript of Philippe Carmona. In this note a large deviation principle for $L_t^i$ was established taking into account the renewal theorem.\medskip
	
	There are two more papers which we would like to mention for continuous space analogue questions. In their analysis of laws of the iterated logarithms for local times of symmetric L\'evy processes, moment generating functions of local times were considered in \cite{MR96}. They exploited the renewal equation (\ref{ren}) where now the transition probabilities need to be replaced by transition kernels. Solving in Laplace domain as we did in part 3. of the proof of Theorem \ref{t1} they transformed back via inverse Laplace transformation to estimate rather delicately the difference 
	\begin{align*}
		 \E^i\big[e^{\gamma L^i_t}\big]- \frac{1}{\hat p^{-1}(1/\gamma)\gamma\int_0^{\infty}se^{-\hat p^{-1}(1/\gamma)s}p_s(i,i)\,ds} e^{\hat p^{-1}(1/\gamma)t}.
	\end{align*}
	Their estimate is uniform in $t$ and $\gamma$ but does not establish convergence as $t$ tends to infinity. Applying our proofs to the renewal equation representation (see the proof of their Lemma 2.6), one obtains the same results for local times of L\'evy processes as we obtained in discrete space.\medskip
	
	Recently, a parabolic Anderson model in $\mathbb{R}$ with L\'evy driver was consider in \cite{FK1} and \cite{FK2}. As their results are based on the same renewal equation (see for instance Equation (2.2) of \cite{FK2} or (4.15) of \cite{FK1}) that we used, one can strengthen their bounds away from the notion of Lyapunov exponents to strong asymptotics with the same expressions for constants and exponential rates as in our discrete setting. This is not surprising as also for their L\'evy process driven version of the parabolic Anderson model the afore mentioned correspondence of second moments and exponential moments of local times of the corresponding L\'evy process holds true.

\end{document}